# Character Sums Over The Prime Numbers
N. A. Carella, March 2012

*Abstract*: A few elementary estimates of a basic character sum over the prime numbers are derived here. These estimates are nontrivial for character sums modulo large $q$. In addition, an omega result for character sums over the primes is also included.

## 1 Introduction

Let $q \geq 1$ be an integer, and let $\chi \neq 1$ be a nonprincipal character modulo $q$. The tasks of determining nontrivial estimates and explicit estimates of the basic character sums

$$\sum_{p \leq x} \chi(p), \quad \text{and} \quad \sum_{p \leq x} \chi(p) \tag{1}$$

over the primes, arithmetic progressions or the integers up to $x \geq 1$ are extensively studied problems in the analytic number theory literature, [KA], [GS], [GK], [GR], [PM], [TR], et alii.

The earliest nontrivial estimate of the basic character sum over the prime numbers, in the literature, seems to be the Vinogradov estimate

$$\sum_{p \leq x} \chi(p) = O(\pi(x) q^{-\delta}), \tag{2}$$

for arbitrarily small real number $\varepsilon > 0$, all sufficiently large number $x > q^{1+\varepsilon}$, and $\delta = \delta(\varepsilon) > 0$, see [KA, p. 156]. This is slightly better than the trivial estimate $\sum_{p \leq x} \chi(p) \leq \pi(x)$, where $\pi(x) = \#\{ p \leq x : p \text{ prime} \}$. Many other estimates of various forms such as $\sum_{p \leq x} \chi(p) = O(x / \log^B x)$, $B > 1$ constant, are available in the literature, see [IK, p. 348], [FR], and similar references. A survey of open problems in exponential and character sums appears in [SK].

The elementary estimates and evaluation in Theorems 1, 2, 3, and 4 for character sums over the prime numbers are derived by elementary methods, and appear to be nontrivial for character sums over the prime numbers modulo large $q$.

## 2 Character Sums Over The Prime Numbers

There are various other related, and more sophisticated results available in the literature, confer the more recent works given in [KA], [FR], et alii. The estimates of character sums over the primes derived here, by elementary methods, are simpler, and similar to the Polya-Vinogradov inequality, and the Paley inequality in forms and the theoretical frameworks.

The proof uses a windowing technique to spread the calculations into two parts: the main component, and an error component. This spreading is intrinsic in the Fourier series of a window function. It closely follows the pattern of the proof of the Polya-Vinogradov inequality using a window function constructed in [PM], and attributed to Landau.

**Theorem 1.** Let $q \geq 3$ be a large integer, and let $\chi \neq 1$ be a nonprincipal character modulo $q$. Then

$$\left| \sum_{p \leq x} \chi(p) \right| = O(q^{1/2+\varepsilon}) \tag{3}$$

for any large number $x \geq 1$ such that $x \leq q^{1-\varepsilon}$, and arbitrarily small $\varepsilon > 0$. In particular,

$$\max_{1 \leq x \leq q^{1-\varepsilon}} \left| \sum_{p \leq x} \chi(p) \right| = O(q^{1/2+\varepsilon}). \tag{4}$$

*Proof*: Utilize the window function $W(t)$, see Lemma 7, to rewrite the character sum over the primes as

$$\begin{aligned}
\sum_{p \leq x} \chi(p) &= \frac{\chi(x)}{2} + \sum_{1 < p < q} \chi(p) W\left(\frac{2\pi p}{q}\right) \\
&= \frac{\chi(x)}{2} + \sum_{1 < p < q} \chi(p) \left( a_0 + \sum_{m \geq 1} \left( a_m \cos \frac{2\pi mp}{q} + b_m \sin \frac{2\pi mp}{q} \right) \right) \\
&= \frac{\chi(x)}{2} + \sum_{1 < p < q} \chi(p) \left( a_0 + \sum_{1 \leq m \leq K} \left( a_m \cos \frac{2\pi mp}{q} + b_m \sin \frac{2\pi mp}{q} \right) \right) \\
&\quad + \sum_{1 \leq p \leq q} \chi(p) \left( \sum_{K < m} \left( a_m \cos \frac{2\pi mp}{q} + b_m \sin \frac{2\pi mp}{q} \right) \right) \\
&= P_K + Q_K + \frac{\chi(x)}{2} + \frac{x}{q} \sum_{1 < p < q} \chi(p),
\end{aligned} \tag{5}$$

where the initial parameters are set to $M = 0$, $N = x$, and $a_0 = x/q$, see Lemma 7. Further, the parity of the character $\chi \ne 1$ classifies the character sum as one of the two types:

$$\rho(\bar{\chi}) \sum_{p \le x} \chi(p) = \begin{cases} \sum_{p \le x} \sum_{1 \le t \le q} \bar{\chi}(t) \cos \dfrac{2\pi pt}{q}, & \text{if } \chi(-1) = 1 \text{ is even,} \\ \sum_{p \le x} \sum_{1 \le t \le q} \bar{\chi}(t) \sin \dfrac{2\pi pt}{q}, & \text{if } \chi(-1) = -1 \text{ is odd,} \end{cases} \qquad (6)$$

where $\rho(\chi) = \sum_{p<q} \chi(p) e^{i2\pi p/q}$ is the complete prime exponential sum, see Theorem 2, and Lemmas 5 and 6. In synopsis, this simplifies the analysis since the sine terms vanish if the character is even; likewise, the cosine terms vanish if the character is odd.

**Case of Even Character $\chi(-1) = 1$.** In this case the expressions $P_K$ and $Q_K$ have the followings upper bounds:

$$\begin{aligned} P_K &= \sum_{1 < p < q} \chi(p) \sum_{m \le K} \left( a_m \cos \frac{2\pi mp}{q} + b_m \sin \frac{2\pi mp}{q} \right) \\ &= \sum_{1 < p < q} \sum_{m \le K} a_m \chi(p) \cos \frac{2\pi mp}{q} \\ &= \rho(\chi) \sum_{m \le K} a_m \bar{\chi}(m) , \end{aligned} \qquad (7)$$

refer to equations (5) and (6). Now, replacing the coefficients $a_m$, see Lemma 7, taking absolute value, and simplifying return

$$\begin{aligned} |P_K| &= \left| \rho(\chi) \sum_{m \le K} a_m \bar{\chi}(m) \right| \\ &= \left| \rho(\chi) \sum_{m \le K} \left( \frac{1}{\pi m} \sin \frac{2\pi mx}{q} \right) \bar{\chi}(m) \right| \\ &\le \frac{q^{1/2}}{\pi} \sum_{m \le K} \frac{1}{m} \\ &\le c_0 q^{1/2} \log K, \end{aligned} \qquad (8)$$

where $|\rho(\chi)| \le \pi(q)^{1/2} \le q^{1/2}$, see Theorem 2, and $c_0 > 0$ is a constant. A slightly different procedure is applied to the expression $Q_K$: replacing the coefficients $a_m$ yields

$$\begin{aligned} Q_K &= \sum_{1 < p < q} \chi(p) \sum_{K < m} \left( a_m \cos \frac{2\pi mp}{q} + b_m \sin \frac{2\pi mp}{q} \right) \\ &= \sum_{1 < p < q} \chi(p) \sum_{K < m} \left( \frac{-1}{\pi m} \cos \frac{2\pi mx}{q} \right) \cos \frac{2\pi mp}{q} \\ &= \frac{-1}{\pi} \sum_{1 < p < q} \chi(p) \sum_{K < m} \frac{1}{m} \left( \sin \frac{2\pi(x-p)m}{q} + \sin \frac{2\pi mx}{q} \right), \end{aligned} \qquad (9)$$

see equation (5), and Lemma 7. Now, taking absolute value, and simplifying return

$$\begin{aligned} |Q_K| &= \left| \frac{1}{\pi} \sum_{1 < p < q} \chi(p) \sum_{K < m} \frac{1}{m} \left( \sin \frac{2\pi(x-p)m}{q} + \sin \frac{2\pi mp}{q} \right) \right| \\ &\leq 2 \left| \sum_{1 < p < q} \chi(p) \sum_{m \leq K} \frac{1}{\pi m} \left( \sin \frac{2\pi mp}{q} \right) \right| \\ &\leq \frac{4}{\pi} \sum_{1 < p < q} \left| \sum_{K < m} \frac{1}{m} \sin \frac{2\pi mp}{q} \right| \\ &\leq \frac{4}{\pi} \sum_{1 \leq n \leq q} \left| \sum_{K < m} \frac{1}{m} \sin \frac{2\pi mn}{q} \right| \\ &\leq c_1 \frac{q}{K+1} \log K, \end{aligned} \qquad (10)$$

where $c_1 > 0$ is a constant, this follows from the estimate of Lemma 8. Rearrange (5) as

$$\sum_{p \leq x} \chi(p) = \frac{q}{q-x}\left( P_K + Q_K + \frac{\chi(x)}{2} \right), \qquad (11)$$

and consider the list of all the estimates:

$$|P_K| \leq c_0 q^{1/2} \log K, \qquad |Q_K| \leq c_1 \frac{q}{K+1} \log K, \qquad \left| \frac{\chi(x)}{2} \right| \leq 1. \qquad (12)$$

Put $K = K(x)$, for example, $K = x^2$ if $x \geq q^{1/2}$, or $K = x^4$ if $x \geq q^{1/4}$, etc. Then, it quickly follows that $q/(q-x) \leq c_2$ for any large number $x \geq 1$ such that $x \leq q^{1-\varepsilon}$, $c_2 > 0$ constant. Lastly, applying the triangle inequality yields

$$\left| \sum_{p \leq x} \chi(p) \right| \leq \left| \frac{q}{q-x}\left( P_K + Q_K + \frac{\chi(x)}{2} \right) \right| \leq c_3 q^{1/2+\varepsilon}, \qquad (13)$$

where $c_3 > 0$ constant.

**Case of Odd Character $\chi(-1) = -1$.** The proof is similar. ∎

For a fixed $q \geq 3$, and nonprincipal character $\chi \neq 1$ modulo $q$, the prime character sum $\sum_{p < x} \chi(p)$ is unbounded as $x \to \infty$, see Theorem 4. But it has an upper bound on the interval $[1, x]$, $x \leq q$. The norm of the complete prime exponential sum $\rho(\chi) = \sum_{p < q} \chi(p) e^{i2\pi p/q}$ is calculated below.

**Theorem 2.** Let $q \geq 3$ be a large integer, and let $\chi \neq 1$ be a character modulo $q$. Then $|\rho(\chi)| \leq \sqrt{q}$.

*Proof*: Let $t \geq 0$ be a parameter. By definition, the norm is given by

$$|\rho_t(\chi)|^2 = \left(\sum_{p<q} \chi(p) e^{i2\pi pt/q}\right)\left(\sum_{r<q} \bar{\chi}(r) e^{-i2\pi rt/q}\right) \qquad (14)$$
$$= \pi(q) + a + \sum_{p,r<q,\, p\neq r} \chi(pr^{-1}) e^{i2\pi(p-r)t/q},$$

where $p$ and $r$ run over the prime numbers up to $q \geq 3$, $\pi(q)$ denotes the number of primes counting function, $\omega(n) = \#\{p \mid n\}$, denotes the number of prime divisors counting function, and

$$a = \begin{cases} \omega(q) - 1 & \text{if } q \text{ is prime,} \\ \omega(q) & \text{if } q \text{ is composite.} \end{cases} \qquad (15)$$

The averages of the norm $|\rho_t(\chi)|^2$ with respect to the variables $t$ and $\chi$ are:

$$E_\chi = \frac{1}{q} \sum_{0<t\leq q} |\rho_t(\chi)|^2 = \pi(q) + a, \qquad \text{and} \qquad E_t = \frac{1}{\varphi(q)} \sum_{\chi \bmod q} |\rho_t(\chi)|^2 = \pi(q) + a \qquad (16)$$

respectively. These information shows that the norm $|\rho_t(\chi)|^2$ of any prime exponential sum $\rho_t(\chi)$ is independent of the variables $t$ and $\chi$. Moreover, the oscillating error term in (14) is a complex number

$$R(t,\chi) = \sum_{p,r<q,\, p\neq r} \chi(pr^{-1}) e^{i2\pi(p-r)t/q} = |R(t,\chi)| e^{i\theta} \qquad (17)$$

where the magnitude $|R(t,\chi)| \geq 0$, and the angle $\theta = \theta(t,\chi)$ are functions of the variables $t$ and $\chi$. Since

$$-|R(t,\chi)| \leq |R(t,\chi)| e^{i\theta} \leq |R(t,\chi)|, \qquad (18)$$

And the norm is nonnegative

$$|\rho_t(\chi)|^2 = \pi(q) + a + |R(t,\chi)| e^{i\theta} \geq 0. \qquad (19)$$

It readily follows that the error term is bounded as $|R(t,\chi)| \leq \pi(q) + a$. Therefore, the norm has the form

$$|\rho_t(\chi)|^2 = \pi(q) + a + R(t,\chi) \leq q, \qquad (20)$$

with the error term $|R(t,\chi)| = O(q/\log q)$. ∎

## 3 Comparisons Of Character Sums Over The Primes And Integers

The two basic character sums (1) are linked via the formula

$$\sum_{p \leq x} \chi(p) = \sum_{n \leq x} \chi(n) \left( -\sum_{d \mid n} \mu(d) \frac{\log d}{\log n} \right) + O(x^{1/2}), \tag{21}$$

so it should not be surprising that these character sums have comparable upper bounds on the interval $[1, x]$. This relationship arises after an application of the Mobius inversion pair

$$\sum_{d \mid n} \Lambda(d) = \log n \quad \text{and} \quad \Lambda(n) = \sum_{d \mid n} \mu(d) \log n/d, \tag{22}$$

and the identity

$$\sum_{p \leq x} \chi(p) = \sum_{n \leq x} \chi(n) \frac{\Lambda(n)}{\log n} + O(x^{1/2}), \tag{23}$$

which decouples the complicated character sum over the primes into a product of two simpler character sums.

Furthermore, the mean square value basic character sum over the primes is given by

$$\sum_{\chi \in S(Q)} \left| \sum_{p \leq Q^\alpha} \chi(p) \right| \leq cQ^{\beta+\varepsilon}, \tag{24}$$

where $S(Q)$ is the set of primitive characters of conductor up to $Q$, $\beta = \max\{1 + \alpha/2, 1/2 + \alpha\}$, $0 < \alpha \leq 1$, and $c > 0$ constant, see [HB] for more general results. Compare this to (15).

The estimate of the basic character sum $\sum_{n \leq x} \chi(n)$ over the integers is given by the Polya-Vinogradov inequality.

**Theorem 3.** Let $q \in \mathbb{N}$ be a large integer, and let $\chi \neq 1$ be a character modulo $q$. Then

$$\left| \sum_{n \leq x} \chi(n) \right| \leq 2q^{1/2} \log q \tag{25}$$

for any real number such that $x > q^{1-\varepsilon}$, $\varepsilon > 0$.

*Proof*: Let $q \not\equiv 2 \bmod 4$, and let $\chi \neq 1$ be a character modulo $q$. The basic character sum over the integers can be written as

$$\sum_{n \leq x} \chi(n) = \sum_{n \leq x} \frac{1}{\tau(\bar{\chi})} \sum_{1 \leq t < q} \bar{\chi}(t) \, e^{i 2\pi nt/q}, \tag{26}$$

see Lemma 6. Rearrange and evaluate this representation into the form

$$\sum_{n \leq x} \chi(n) = \frac{1}{\tau(\overline{\chi})} \sum_{1 \leq t < q} \overline{\chi}(t) \left( \frac{e^{i2\pi(x+1)t/q} - 1}{e^{i2\pi t/q} - 1} \right), \tag{27}$$

where $|\tau(\chi)| = \sqrt{q}$, or $|\tau(\chi)| = \sqrt{2q}$. An upper estimate is

$$\left| \sum_{n \leq x} \chi(n) \right| \leq \frac{1}{\sqrt{q}} \sum_{1 \leq t < q} \frac{1}{|\sin \pi t/q|} \tag{28}$$
$$\leq cq^{1/2} \log q,$$

where $|\sin \pi t/q| \geq \pi t/q$ for $1 \leq t < q$ was utilized, and $c > 0$ is a constant. ∎

The exceptional value $q \equiv 2 \bmod 4$ stems from the vanishing of the Gaussian sum $\sum_{0 < t < p} \chi(t) e^{i2\pi t/q} = 0$ whenever $\chi \neq 1$ is a quadratic character modulo $q \equiv 2 \bmod 4$, see Lemma 6. The numbers $q \equiv 2 \bmod 4$ seem to be exceptional values of the Polya-Vinogradov inequality $\sum_{n \leq x} \chi(n) \leq cq^{1/2} \log q$, $c > 0$ constant. But this seems to be irrelevant since there is different proof of the Polya-Vinogradov inequality, refer to [DW].

The Polya-Vinogradov estimate is nontrivial for any real number such that $x > q^{1-\varepsilon}$, $\varepsilon > 0$. The estimates of character sums over short intervals of length $x < q^{1-\varepsilon}$ are more delicate and complex. The analysis for the short range of values $x < q^{1-\varepsilon}$ are given in [GS], [GK], [GR], [BR], LZ] and similar literature. The nontrivial estimates of $\sum_{n \leq x, n \equiv a \bmod q} \chi(n)$ over arithmetic progressions $\{ p = qn + a : \gcd(a, q) = 1, \text{ and } n \geq 1 \}$ are given in [FR].

## 4 Omega Result For Character Sums Over The Prime Numbers

The Polya-Vinogradov inequality is about the best possible since the Paley inequality

$$\sum_{n \leq x} \chi(n) \geq cq^{1/2} \log \log q \tag{29}$$

holds for some characters $\chi \neq 1$, and infinitely many primes $q$, and $c > 0$ constant, confer [MV, p. 312] for a proof. A similar omega result can be obtained for the basic character sums over the prime numbers.

Let $q$ be a large integer, and let $g$ be a primitive root modulo $q$. For $q \neq 2^m$, a multiplicative character has the form $\chi_k(n) = e^{i2\pi k \, \mathrm{Ind}_g n/\varphi(q)}$, where $0 \leq k < q$, and $\mathrm{Ind}_g(n)$ is the discrete logarithm modulo $q$, see Section 6. The character $\chi$ is called even if $\chi(-1) = 1$, otherwise, it is odd and $\chi(-1) = -1$. For example, every character such that $\gcd(2, k) = 1$ is odd.

The cancelation mechanism of some Dirichlet characters springs from some unique pairwise partitions of the multiplicative group of the integers $\mathbb{Z}_q$ modulo $q$. Let the congruence $p \equiv g^v \bmod q$, $0 \leq v < \varphi(q) - 1$, specifies the $\varphi(q)$ equivalent classes of primes modulo $q$. Some unique pairwise partitions of the multiplicative group modulo $q$ are as follows:

(i) For odd character $\chi(-1) = -1$ modulo $q = 2m + 1$, the pairwise partition is

$$\left(g^v, g^{v+\phi(q)/2}\right), \quad \text{for } 0 \leq v < \varphi(q)/2, \tag{30}$$

(ii) For quartic character $\chi(\pm i) = -1$ modulo $q = 4m + 1$, the pairwise partition is

$$\left(g^v, g^{v+\phi(q)/4}\right), \left(g^{v+1+\phi(q)/2}, g^{v+1+3\phi(q)/4}\right), \quad \text{for } 0 \leq v < \varphi(q)/4,$$

and similar pairwise partitions for some other $q = 2^a m + 1$, and $a > 2$, mutatis mutandis. In addition, the prime counting function on the arithmetic progression $\{ p = qn + a : \gcd(a, q) = 1, \text{ and } n \geq 1 \}$ is defined by

$$\pi(x, a, q) = \#\{ p \leq x : p \equiv a \bmod q \} = \varphi(q)^{-1} li(x) + E(x, a, q), \tag{31}$$

where $li(x) = \int_1^x (\log t)^{-1} dt$ is the logarithmic integral, and $E(x, a, q)$ is the prime number theorem error term.

**Theorem 4.** Let $q \in \mathbb{N}$ be a fixed integer, and let $\chi \neq 1$ be an odd character modulo $q$. Then

$$\sum_{p \leq x} \chi(p) = \Omega_\pm \left( \frac{x^{1/2} \log\log\log x}{\log x} \right) \tag{32}$$

for some sufficiently large real number $x \geq q^{1-\varepsilon}$, and $\varepsilon > 0$ arbitrarily small. In particular,

$$\sum_{p \leq x} \chi(p) = \Omega_\pm \left( q^{1/2} \log\log\log q / \log q \right) \tag{33}$$

for infinitely many large $q$ and odd character.

Proof: Let $g$ be a primitive root modulo $q$, and choose an odd character $\chi(g^{\varphi(q)/2}) = \chi(-1) = -1$. Since $\chi(g^{v+\varphi(q)/2}) = \chi(g^v)\chi(g^{\varphi(q)/2}) = -\chi(g^v)$, the partition pairing $g^v$, $g^{v+\varphi(q)/2}$ of residue classes, see (10?), yields

$$\pi(x, g^v, q)\chi(g^v) + \pi(x, g^{v+\varphi(q)/2}, q)\chi(g^{v+\varphi(q)/2}) = \pi(x, g^v, q)\chi(g^v) - \pi(x, g^{v+\varphi(q)/2}, q)\chi(g^v)$$
$$= \left( E(x, g^v, q) - E(x, g^{v+\varphi(q)/2}, q) \right) \chi(g^v). \tag{34}$$

Thus, it follows that the main terms $\pi(x, g^v, q) \sim \varphi(q)^{-1} li(x)$ and $\pi(x, g^{v+\varphi(q)/2}, q) \sim \varphi(q)^{-1} li(x)$ cancel pairwise, $0 \leq v < \varphi(q)/2$. As a consequence, the basic character sum over the primes collapses to

$$\sum_{p \leq x} \chi(p) = \sum_{0 \leq v < \varphi(q)-1} \sum_{p \leq x, \, p \equiv g^v \bmod q} \chi(p)$$
$$= \sum_{0 \leq v < \varphi(q)-1} \chi(g^v) \pi(x, g^v, q) \tag{35}$$
$$= \sum_{0 \leq v < \varphi(q)/2} \chi(g^v) E(x, g^v, q),$$

which is a complex linear combination of the error terms $E(x, g^v, q)$ of the $\varphi(q)$ equivalent classes of primes modulo $q$. Now, the result follows from the Littlewood form of the prime number theorem

$$\pi(x, a, q) = \frac{li(x)}{\varphi(q)} + \Omega_\pm\left(\frac{x^{1/2} \log\log\log x}{\varphi(q)\log x}\right), \qquad (36)$$

on the arithmetic progression $\{ p = qn + a : \gcd(a, q) = 1, \text{ and } n \geq 1 \}$. Refer to [MV, p. 479], [IV, p. 51], and similar literature. ∎

A few examples were computed to demonstrate the concept, and the challenges faced in estimating small character sums over the primes. Here, the congruence $p \equiv g^v \bmod q$, $0 \leq v < \varphi(q) - 1$, specifies the $\varphi(q)$ equivalent classes of primes modulo $q$, and the corresponding odd character $\chi(g^{\varphi(q)/2}) = \chi(-1) = -1$.

1. $\sum_{p \leq x} \chi(p) = \pi(x, 1, 3) - \pi(x, 2, 3) = \Omega_\pm(x^{1/2} \log\log\log x / \log x)$, for $q = 3$, $\chi_3(2^v) = (-1)^v$, $g^v \equiv 2^v \bmod 3$.

2. $\sum_{p \leq x} \chi(p) = \pi(x, 1, 4) - \pi(x, 3, 4) = \Omega_\pm(x^{1/2} \log\log\log x / \log x)$, for $q = 4$, $\chi_4(3^v) = (-1)^v$, $g^v \equiv 3^v \bmod 4$.

3. $\sum_{p \leq x} \chi(p) = \pi(x, 1, 5) + i\pi(x, 2, 5) - i\pi(x, 3, 5) - \pi(x, 4, 5) = \Omega_\pm(x^{1/2} \log\log\log x / \log x)$, (37)

for $q = 5$, $\chi_5(2^v) = i^v$, $g^v \equiv 2^v \bmod 5$.

Similar, and other related advanced topics in comparative number theory, are discussed in [RN, p. 275], [GM], [MJ], [RS], [SC], et alii.

***Theorem 5.*** Let $q \in \mathbb{N}$ be an integer, and let $\chi \neq 1$ be a character modulo $q$. Then

$$\sum_{p \leq x} \chi(p) = O\left(xe^{-c(\log x)^{1/2}}\right) \qquad (38)$$

for some sufficiently large real number $x \geq x_0$, and $c > 0$ constant.

Proof: Since $\sum_{1 \leq a \leq \varphi(q)} \chi(a) = 0$, the basic character sum over the primes collapses to

$$\begin{aligned}\sum_{p \leq x} \chi(p) &= \sum_{1 \leq a \leq \varphi(q),\ p \leq x, p \equiv a \bmod q} \sum \chi(p) \\ &= \sum_{1 \leq a \leq \varphi(q)} \chi(a)\pi(x, a, q) \qquad (39) \\ &= \sum_{1 \leq a \leq \varphi(q)} \chi(a)E(x, a, q),\end{aligned}$$

which is a complex linear combination of the error terms $E(x, a, q)$ of the $\varphi(q)$ equivalent classes of primes modulo $q$. Now, the result follows from the delaVallee Poussin form of the prime number theorem

$$\pi(x, a, q) = \frac{li(x)}{\varphi(q)} + O\left(xe^{-c(\log x)^{1/2}}/\varphi(q)\right), \tag{40}$$

on the arithmetic progression $\{ p = qn + a : \gcd(a, q) = 1, \text{ and } n \geq 1 \}$. Refer to [MV, p. 479], [IV, p. 51], and similar literature. ∎

## 5. Elementary Foundation
This Section serves as a reference of some of the concepts used to complete the proof of the main results in previous Section.

### 5.1 A Few Arithmetic Functions
Let $n = p_1^{v_1} \cdot p_2^{v_2} \cdots p_r^{v_r}$, let $\omega(n) = \#\{ p \mid n \}$ be the number of prime divisors counting function, and let $\Omega(n) = v_1 + v_2 + \cdots + v_r$. The Mobius function defined by

$$\mu(n) = \begin{cases} (-1)^{\omega(n)} & \text{if } \omega(n) = \Omega(n), \\ 0 & \text{if } \omega(n) \neq \Omega(n), \end{cases} \tag{41}$$

the vonMangold function defined by

$$\Lambda(n) = \begin{cases} \log p & \text{if } n = p^k, \\ 0 & \text{if } n \neq p^k, \end{cases} \tag{42}$$

where $p^k$ is a prime power. And let the Euler function be defined by

$$\varphi(q) = \prod_{p \mid q} (1 - 1/p), \tag{43}$$

where $p$ ranges over the prime divisors of $q$. A few other number theoretical functions are also used throughout the paper.

### 5.2 Multiplicative Characters
A multiplicative character $\chi$ is a periodic, complex valued and completely multiplicative function $\chi : \mathbb{Z} \to \mathbb{C}$ on the integers. For each $q \in \mathbb{N}$, the set of characters $\hat{G} = \{ 1 = \chi_0, \chi_1, \ldots, \chi_{\varphi(q)-1} \}$ is a group of order $\varphi(q)$.

**Lemma 6.** If a function $f : \mathbb{Z} \to \mathbb{C}$ satisfies $f(n) \equiv 0 \mod q$ for $\gcd(n, q) \neq 1$, is periodic, and completely multiplicative, then $f(n) = \chi(n)$ is a character modulo $q$.

**Properties of Nontrivial Characters**
(i) $\chi(1) = 1$ and $\chi(s) \neq 1$,     for $\gcd(s, q) > 1$,
(ii) $\chi(st) = \chi(s)\chi(t)$,     multiplicative,
(iii) $\chi(s + k) = \chi(s)$,     periodic of period $k \geq 1$,
(iv) $|\chi(s)| = 1$,     a point in unit circle.

There are several forms of the multiplicative characters depending on the decomposition of the integer $q$.

For $q \neq 2^m$, $m \geq 1$. There are $\varphi(q)$ characters in $\hat{G} = \{1 = \chi_0, \chi_1, \ldots, \chi_{\varphi(q)-1}\}$, the principal character $\chi = 1$ and nonprincipal characters $\chi \neq 1$ are defined by

$$\chi_0(n) = \begin{cases} 1 & \text{if } \gcd(n,q) = 1, \\ 0 & \text{if } \gcd(n,q) \neq 1, \end{cases} \tag{44}$$

and

$$\chi_k(n) = \begin{cases} e^{i2\pi k \log(n)/\varphi(q)} & \text{if } \gcd(n,q) = 1, \\ 0 & \text{if } \gcd(n,q) \neq 1, \end{cases} \tag{45}$$

respectively. The notation $\log_g(n) = \text{Ind}_g(n)$ denotes the discrete logarithm with respect to some primitive root $g \bmod q$.

For $q = 2^v$, $v \geq 1$. A character is realized by one of the three forms described below.
Case $v = 1$, there is a single character in $\hat{G} = \{1 = \chi_0\}$, and it is defined by

$$\chi(n) = \begin{cases} 1 & \text{if } n \equiv 1 \bmod 2, \\ 0 & \text{otherwise.} \end{cases} \tag{46}$$

Case $v = 2$. There are two characters in $\hat{G} = \{1 = \chi_0, \chi\}$, the nontrivial character is defined by

$$\chi(n) = \begin{cases} 1 & \text{if } n \equiv 1 \bmod 4, \\ -1 & \text{if } n \equiv 3 \bmod 4, \\ 0 & \text{if } n \equiv 0, 2 \bmod 4. \end{cases} \tag{47}$$

Case $v > 2$. There are $\varphi(q)$ characters in $\hat{G} = \{1 = \chi_0, \chi_1, \ldots, \chi_{\varphi(q)-1}\}$, and a nontrivial character is defined by

$$\chi_s(n) = \begin{cases} (-1)^\delta e^{i2\pi st/\varphi(2^{v-3})} & \text{if } n \equiv 1 \bmod 2, \\ 0 & \text{if } n \equiv 0 \bmod 2, \end{cases} \tag{48}$$

for some $0 \leq s, t < \varphi(q)$. The integer $n$ is represented as $n \equiv (-1)^\delta 5^t \bmod 2^v$ in the multiplicative group of units $\{-1, 1\} \times \{5^t : 0 \leq t < \varphi(2^{v-2})\}$ of $\mathbb{Z}_q$, where $\delta = 0$ if $n \equiv 1 \bmod 4$ or $\delta = 1$ if $n \equiv 3 \bmod 4$. This is due to the fact that this multiplicative group is not cyclic.

A character $\chi$ is *even* if $\chi(t) = \chi(-t)$, otherwise $\chi(t) = -\chi(-t)$, and the character is *odd*. The binary variable $\delta_\chi = 0, 1$ tracks the even odd condition, specifically, $\chi(-n) = (-1)^{\delta_\chi} \chi(n)$. A character $\chi$ is *primitive* if no proper subgroup of the group $\hat{G} = \{1 = \chi_0, \chi_1, \ldots, \chi_{\varphi(q)-1}\}$ contains it. Under this condition the conductor of a character is the integer $f_\chi = q$.

**Lemma 7.** (Orthogonal relations) For $a \geq 1$, and the set of characters modulo $q$, the followings hold.

(i) $\displaystyle\sum_{1 \leq n \leq q} \chi(n) = \begin{cases} \varphi(q) & \text{if } \chi = \chi_0, \\ 0 & \text{if } \chi \neq \chi_0. \end{cases}$  (ii) $\displaystyle\sum_{\chi \bmod q} \chi(n) = \begin{cases} \varphi(q) & \text{if } n \equiv 0 \bmod \varphi(q), \\ 0 & \text{otherwise}. \end{cases}$ (49)

A Gauss sum is defined by the exponential sum

$$\tau_a(\chi) = \sum_{n=1}^{q} \chi(n) e^{i2\pi a n/q}, \quad 0 \leq a < q. \tag{50}$$

**Lemma 8.** Let $\chi \neq \chi_0$ be a nontrivial character modulo $q$, and let $\tau(\chi) = \tau_1(\chi)$. Then
(i) $\tau_a(\chi) = \chi(a)\tau(\overline{\chi})$,  (ii) $\tau_a(\chi)\tau_a(\overline{\chi}) = q$,

(iii) $\tau_a(\chi) = \begin{cases} (1+i)q^{1/2} & \text{if } q \equiv 0 \bmod 4, \\ q^{1/2} & \text{if } q \equiv 1 \bmod 4, \\ 0 & \text{if } q \equiv 2 \bmod 4, \\ iq^{1/2} & \text{if } q \equiv 3 \bmod 4. \end{cases}$  (51)

Let $q \not\equiv 2 \bmod 4$, and let $\chi \neq 1$ be a nonprincipal character modulo $q$. The Fourier transform

$$\chi(n) = \frac{1}{\tau(\overline{\chi})} \sum_{1 \leq t < q} \overline{\chi}(t) e^{-i2\pi n t/q} \tag{52}$$

is a complex-valued function that interpolates the Dirichlet character $\chi$ form $\mathbb{N}$ to $\mathbb{C}$. For $q \equiv 2 \bmod 4$, it is undefined. The interpolation formula has a simpler form identified by the parity of the character:

$$\tau(\overline{\chi})\chi(n) = \begin{cases} \displaystyle\sum_{1 \leq t < q} \overline{\chi}(t) \cos \frac{2\pi n t}{q}, & \text{if } \chi(-1) = 1 \text{ is even}, \\ \displaystyle\sum_{1 \leq t < q} \overline{\chi}(t) \sin \frac{2\pi n t}{q}, & \text{if } \chi(-1) = -1 \text{ is odd}. \end{cases} \tag{53}$$

### 5.3 Window Function, and Estimate of Trigonometric Series
Windowing schemes are widely used in signal analysis and number theory. For a pair of integers $M < N$, defines the window function

$$W(x) = \begin{cases} 1/2, & \text{if } x = 2\pi M/q \text{ or } x = 2\pi N/q, \\ 1, & \text{if } 2\pi M/q < x < 2\pi N/q, \\ 0, & \text{otherwise}, \end{cases} \qquad (54)$$

on the interval $[0, 2\pi)$.

**Lemma 9.** The Fourier series of the window function $W(x)$ is given by

$$W(x) = a_0 + \sum_{m \geq 1} (a_m \cos mx + b_m \sin mx), \qquad (55)$$

where the first coefficient is $a_0 = (N - M)/q$, and for $m \geq 1$, the coefficients are given by

$$a_m = \frac{1}{\pi m}\left(\sin\frac{2\pi Nm}{q} - \sin\frac{2\pi Mm}{q}\right), \quad \text{and} \quad b_m = \frac{-1}{\pi m}\left(\cos\frac{2\pi Nm}{q} - \cos\frac{2\pi Mm}{q}\right). \qquad (56)$$

The pair of functions (48) and (49) is a Fourier pairs. Other well known Fourier pairs are the followings:

$$f(t) = \begin{cases} 0, & \text{if } |t| > 1, \\ 1, & \text{if } |t| < 1, \\ 1/2, & \text{if } |t| = 1, \end{cases} \quad \text{and} \quad \hat{f}(s) = \frac{\sin 2\pi s}{\pi s}, \qquad (57)$$

and

$$g(t) = \max\{1 - |t|, 0\} \quad \text{and} \quad \hat{g}(s) = \left(\frac{\sin \pi s}{\pi s}\right)^2, \qquad (58)$$

for $t \in \mathbb{R}$, and $s \in \mathbb{C}$. Other versions and proofs of Lemmas 7 and 8 are discussed in [PM, p. 4].

**Lemma 10.** Let $q \geq 3$ be an integer, and let $K > 1$ be a large number. Then

$$\sum_{1 \leq n \leq q} \left| \sum_{K < m} \frac{1}{m} \sin(mn/q) \right| \leq c \frac{q}{K+1} \log K, \qquad (59)$$

where $c > 0$ is a constant.